\documentclass[11pt]{article}
\usepackage{amsmath}
\usepackage{amssymb}
\usepackage{amscd}
\usepackage{amsfonts}

\def\reel{\hbox{{\rm R}\kern-1em\hbox{{\rm I} }}}
\def\relatif{\ \hbox{{\rm Z}\kern-.4em\hbox{\rm Z}}}
\def\nat{\hbox{{\rm N}\kern-1em\hbox{{\rm I} } }}
\def\comp{\hbox{{\rm C}\kern-.55em\hbox{{\rm I} } }}
\def\smallcomp{\hbox{\fiverm C}\kern-.35em{\hbox{\fiverm I}}}

\def\fudge{\mathchoice{}{}{\mkern.5mu}{\mkern.8mu}}
\def\bbc#1#2{{\rm \mkern#2mu\vbar\mkern-#2mu#1}}
\def\bbb#1{{\rm I\mkern-3.5mu #1}} \def\bba#1#2{{\rm #1\mkern-#2mu\fudge
#1}}
\def\bb#1{{\count4=`#1 \advance\count4by-64 \ifcase\count4\or\bba
A{11.5}\or \bbb B\or\bbc C{5}\or\bbb D\or\bbb E\or\bbb F
\or\bbc G{5}\or\bbb H\or \bbb I\or\bbc J{3}\or\bbb K\or\bbb
L \or\bbb M\or\bbb N\or\bbc O{5} \or \bbb P\or\bbc
Q{5}\orrrr \bb R\or\bbc S{4.2}\or\bba T{10.5}\or\bbc
U{5}\or    \bba V{12}\or\bba W{16.5}\or\bba X{11}\or\bba
Y{11.7}\or\bba Z{7.5}\fi}}
\def\rat{\hbox{{\rm Q}\kern-.70em\hbox{{\rm I} } }}
\def \R {\bbb R}

\def \P {\bbb P}
\def \E {\bbb E}
\newcommand{\refm [1]} {(\ref{#1})}
\def\calF{{\cal F}}

\newcounter{corcountrer}
\newcounter{theoremcounter}
\newcounter{lemmacounter}
\newcounter{problemcounter}
\newcounter{remarkcounter}
\newcounter{propositioncounter}
\newcounter{definitioncounter}
\newtheorem{cor}[corcountrer]{Corollary}
\newtheorem{thm}[theoremcounter]{Theorem}
\newtheorem{lemma}[lemmacounter]{Lemma}
\newtheorem{prob}[problemcounter]{Problem}
\newtheorem{remark}[remarkcounter]{Remark}
\newtheorem{proposition}[propositioncounter]{Proposition}
\newtheorem{definition}[definitioncounter]{Definition}

\newcommand{\be}{\begin{equation}}
\newcommand{\ber}{\begin{eqnarray}}
\newcommand{\ql}{\mathbb{R}^l}
\newcommand{\nin}{\noindent}

\def\qed{\hfill \vrule height1.3ex width1.2ex depth-0.1ex}
\def\bbb#1{{\rm I\mkern-3.5mu #1}} \def\bba#1#2{{\rm #1\mkern-#2mu\fudge
#1}}

\newcommand{\la}{\label}
\newcommand{\en}{\end{equation}}

\title{Asymptotics of counts of small components in random  structures
and  models of coagulation-fragmentation. }
\author{{\bf Boris L. Granovsky}
\thanks{E-mail: mar18aa@techunix.technion.ac.il} \\
Department of Mathematics, Technion-Israel Institute of Technology,\\
Haifa, 32000,Israel.}

\begin{document}
\maketitle \vskip 5cm
\nin American Mathematical Society 1991 subject classifications.

\nin Primary-60C05;60K35, secondary-05A16,82B05, 11M45,

\nin Keywords and phrases: Multiplicative measures on the set of
partitions, Random structures, Coagulation-Fragmentation
processes, Schur's lemma, Models of ideal gas.



\begin{center}
{\bf Abstract}
\end{center}
 \nin We establish necessary and sufficient conditions for
the convergence (in the sense of finite dimensional distributions) of
 multiplicative measures on the set of partitions.
 The multiplicative
measures depict  distributions of component spectra of random
structures and also the equilibria of classic models of
statistical mechanics and
 stochastic processes of coagulation--fragmentation.
 We show that
  the convergence of multiplicative measures is equivalent to the asymptotic
  independence of counts of
  components of fixed sizes in random structures.
 We then apply Schur's tauberian lemma
and some  results from additive number theory and enumerative
combinatorics in order to derive plausible sufficient conditions of convergence. Our results demonstrate
that the common belief, that counts  of components of fixed sizes
in random structures  become
independent as the number of particles goes to infinity, is not
true in general.

\newpage

\nin \section{Introduction: Probabilistic setting and its applications}
 \setcounter{equation}{0}
We start from the following formalism.

\nin  Let $\{Z_j,\ j\ge 1\}$ be a sequence of independent integer
valued random variables that induces a sequence of random vectors
$\{{\bf K}^{(n)}=(K_1^{(n)},\ldots,K_n^{(n)}),\ n\ge 1\}$ given by

\be {\cal L}({\bf K}^{(n)})= {\cal L}(Z_1,\ldots, Z_n \vert
\sum_{j=1}^njZ_j=n), \quad n=1,2,\ldots. \la{1}
\end{equation}

\nin It follows from \refm[1] that ${\bf K}^{(n)}\in \Omega_n, \
n\ge 1,$ where \be \Omega_n=\{\eta=(k_1,\ldots, k_n):\sum_{j=1}^n
jk_j=n \}
\end{equation}
\nin is the set of all partitions $\eta$ of  an integer $n.$
 \nin In probabilistic combinatorics, \refm[1] is called
the conditioning relation (see\cite{ABT}), while the sequence of
vectors $\{{\bf K}^{(n)},\ n\ge 1\},$ is called the counting
process.

\nin Next, denote by $\mu_n$ the probability measure on $\Omega_n$
induced by the conditioning relation \refm[1]:$$\mu_n(\eta):=\P
({\bf K}^{(n)}=\eta),\quad \eta\in \Omega_n,\quad n\ge 1$$ and let
\be a_k^{(j)}=\P (Z_j=k), \quad k\ge 0, \quad j\ge 1. \la{2}
\end{equation}
\nin We then have \be \mu_n(\eta)=c_n^{-1} \prod_{j=1}^n
a_{k_j}^{(j)}, \quad \eta=(k_1,\ldots, k_n)\in \Omega_n, \quad
n\ge 1, \nonumber
\end{equation}

\nin where \be c_n=\P(\sum_{j=1}^n jZ_j=n)=  \sum_{\eta\in
\Omega_n} \prod_{j=1}^n a_{k_j}^{(j)}, \quad \eta=(k_1,\ldots,
k_n)\la{3}
\end{equation}
\nin is the  partition function for the measure $\mu_n.$ Taking into account \refm[1], \refm[3] we will
assume throughout this paper that the probabilities $a_{k}^{(j)},\
k\ge 0, \ j\ge 1$ \nin are such that  $c_n>0,\ n\ge 1$. Vershik (\cite{V1})
suggested that  the class of measures (\ref{3}) be called
multiplicative, while Pitman (\cite{P}) calls them Gibbsian (see also \cite{GK}). Observe that the multiplicative form \refm[3] of
the sequence of measures $\mu_n$ is implied by the fact that the
random variables $Z_j, \ j\ge 1$ in \refm[1] do not depend on $n.$

 \nin It is clear that the sequence of measures $\{\mu_n, \ n\ge
1\}$ induced by \refm[1]
 is uniquely defined by
the array of probabilities $\{a_k^{(j)},\ j\ge 1,\ k\ge 0\}.$
However, this correspondence is not a bijection. In fact, the
``tilting" transformation (see \cite{ABT}) of the probabilities:
\be a_k^{(j)}(\rho)= \frac{\rho^{jk}a_k^{(j)}}{S^{(j)}(\rho)},\
k\ge 0,\ j\ge 1,\label{ti}\end{equation} where $\rho> 0$ and
$S^{(j)}(\rho)$ is the normalizing constant, does not change the
sequence of measures $\{\mu_n\,\ n\ge 1\}$. But this transformation does
affect the generic partition function $c_n$ leading to the tilted partition function $c_n(\rho):$ \be c_n(\rho)=\frac{c_n
\rho^n}{ \prod_{j=1}^nS^{(j)}(\rho)}, \quad n\ge
1.\label{zz}\end{equation} Note
 that the tilting  is defined for all finite $\rho> 0$ such that
\be S^{(j)}(\rho)=\sum_{k=0}^\infty\rho^{jk}a_k^{(j)}<\infty, \
j\ge 1.\la{sj}\end{equation}

 \vskip .5cm It is a remarkable fact  that
the representation \refm[1] provides a mathematical formalism for
a variety of models in seemingly unrelated contexts. Let us
briefly describe    four main  fields of application of this
setting.

$\bullet$ {\bf Decomposable combinatorial structures} (for more
details see \cite{ABT},\cite{bol},\cite{Kol} and references therein). The size of such a
structure is defined to be the number of elements in it. A
decomposable structure of size $n$ is a union of indecomposable
 components (=components), so that the counts $k_1,\ldots ,k_n$ of
components of sizes $1,\ldots,n $ respectively, form an integer
partition of $n.$ It is  assumed that each  component of size $j$
belongs to one of  $m_j$ types. The
 three classes
of decomposable structures: assemblies, multisets and selections,
encompass the whole  universe of classical combinatorial objects.
Assemblies are structures composed of labeled elements. The class
of assemblies  includes permutations  decomposed into cycles
($m_j=(j-1)!$), forests composed of rooted trees with labeled
vertices
  ($m_j=j^{j-1}$), graphs composed of connected subgraphs with
labeled vertices ($m_j\sim 2^{{j\choose 2}}$), etc. We note that for the last model $m_j$ appears to be asymptotically equal to the  total number  $2^{{j\choose 2}}$ of graphs on $j$ vertices. This follows from the remarkable fact that a random graph on $n$ vertices is connected with probability $1$, as $n\to \infty$. Multisets are
formed from unlabeled elements. Examples of multisets are integer
partitions ($m_j=1$), planar partitions ($m_j=j$) and mapping
patterns ($m_j\sim \frac{\rho^{-j}}{2j},\
 \rho=0.3383$). Regarding the last example, recall that a mapping  from the set $[1,n]$ to itself is a digraph with edges $(i,f(i), \ i=1,\ldots,n)$ decomposed into connected subgraphs of the underlying undirected graph. Mapping patterns are obtained from the above structure by removing labels, so that only the topology of  the graph matters.     Finally, selections are defined as multisets with distinct
components, which means that all component counts $k_j, \
j=1,\ldots,n$ are either $0$ or $1.$ A typical example of a
selection is an integer partition into distinct parts ($m_j=1$).

The basic problem in enumerative combinatorics is to find the
asymptotics, as $n$ goes to infinity, of the number of a certain
class of structures  of size $n$,  with a component spectrum
$(k_1,\ldots,k_n)$ in a given  subset of $\Omega_n$. As a part of
this problem, the asymptotics of the total number of given
structures of size $n$ is of special interest.

 The  starting point of the probabilistic method considered is the definition of
  a random structure of size $n,$ which is a random element
$\Pi_n$ distributed uniformly on the finite set of all structures
considered, with  size  $n$. Next is defined the induced random
component spectrum ${\bf K}^{(n)},$ also called the counting
process:
$${\bf K}^{(n)}= {\bf
K}^{(n)}(\Pi_n)=(K_1^{(n)},\ldots,K_n^{(n)}),\ n\ge 1,\quad {\bf
K}^{(n)}(\Pi_n)\in\Omega_n ,$$ where the random variable
$K_j^{(n)}$ represents the number of components of size $j$ in
$\Pi_n$.

It turns out that
 the representation \refm[1] of the distribution of ${\bf K}^{(n)}$
is valid for the aforementioned three classes of  combinatorial
structures. Namely, assemblies, multisets and selections are
induced respectively, by the following
 three types of random variables $Z_j, \ j\ge 1:$
Poisson $(Po(a_j), \ a_j={m_j\over j!})$, Negative Binomial
$(NB(p^j,m_j))$ and Binomial $(Bi(\frac{p^j}{1+p^j},m_j)), \
0<p<1$.

$\bullet$ {\bf  Models of ideal gas} (for references see
\cite{V1},\cite{Grei}, Ch.12, \cite{Sal}).

 In classical statistical mechanics, an ideal gas is a collection of
perfectly elastic particles (atoms or molecules) which collide but
otherwise do not interact with each other. It is assumed that the
total internal energy $E$ of a gas is the sum of the microscopic
energies of  random motions of individual particles and that $E$
is partitioned between the particles, so that $k_j,$ called an
occupation number, is the number of particles with the energy
level $j,$ having a prescribed  weight $m_j$ that regularly varies
with $j$.
 The following three basic models (=statistics) of ideal gas are common.
 \\ Maxwell-Boltzmann($MB$) (=labeled
particles), Bose-Einstein($BE$) (=indistinguishable particles),
Fermi-Dirac($FD$) (=indistinguishable particles, such that no more
than one particle may have a given energy level). In accordance
with the setting  for combinatorial structures, $MB$, $BE$ and $FD$
models conform to assemblies, multisets and selections,
respectively.  The
probability distribution of the energy states $\eta$ which varies
from model to model, is defined by a measure on the state space
$\Omega_n$.
By laws of statistical mechanics, these measures are
forced to be of the multiplicative form (\ref{3}), with the
numbers $a_k^{(j)}$ defining the type of a model of the ideal gas
considered.\\
A substantial difference of the model of  ideal gas, treated as a
quantum system is that a particle of the $d$-dimensional gas  is
viewed
 as a lattice point
  ${\bf q}\in{\cal Z}_d$ and the energy levels $\epsilon_{\bf q}$,
  called energy eigenvalues,
   are of the following special
  form:
  \be\epsilon_{\Vert \bf{q}\Vert^2}=c\Vert {\bf q}\Vert^2,
  \ {\bf q}=(l_1,\ldots,l_d)\in {\cal Z}^d, \quad \Vert
  {\bf q}\Vert^2=\sum_{s=1}^d l_s^2,
  \la{eps}\end{equation}
where  $c>0$ is a known constant that does not depend on $\bf{q}.$
Consequently, the state of the quantum system is determined by a
weighted  partition $\eta=(k_1,\ldots,k_n)$ of an integer $n=E:$
$E=\sum_{j\ge 1}jk_j.$ By  \refm[eps], to each energy
level $j$ is naturally prescribed a "weight" $r_d(j)$ which is the
number of representations  of the natural number $j$ as the sum of
$d$ integer squares. In other words, $r_d(j)$ is the number of
distinct lattice positions ${\bf q}=(l_1,\ldots,l_d)\in {\cal
Z}^d$ of a particle on a sphere of radius $\sqrt{j} ,$ i.e.
 with the energy level $j.$ It is known from  number theory (see
e.g. \cite{g}) that for $d\ge 5,$ $$C_1j^{\frac{d-2}{2}}\le
r_d(j)\le C_2j^{\frac{d-2}{2}},\ j\ge 1, $$ where $C_1,C_2$ are
positive constants depending on $d,$ and  that for $d=2,3,4$ the
functions $r_d(j)$ oscillate wildly (in $j$), while, obviously,
\be r_1(j)=
  \begin{cases}
    2 , & \text{if  $j$ is a square} \\
    0 , & \text{otherwise}.
  \end{cases}
  \end{equation}

Employing  known  properties of   $r_d(j),$ an important fact was
proven in \cite{V2} that
 for the sake of asymptotic analysis, it is possible to treat
 the $d$- dimensional quantum models as the classic $BE$ and $FD$ ones
 with parameters $m_j=cj^\beta,$ where $ c>0,$
and $\beta=\frac{d-2}{2}, $ if \ $d\ge 2.$

$\bullet$ {\bf Coagulation--fragmentation processes on the set of
integer partitions} (see \cite{DGG},\cite{ABT}).

 \nin  We will show  that a multiplicative measure
$\mu_n$ can be viewed as an equilibrium of a classic
coagulation--fragmentation process $(CFP)$ which is a
time-continuous Markov chain on $\Omega_ n$, defined as follows. A
state $\eta=(k_1,\ldots,k_n)\in\Omega_n$ of a  $CFP$ depicts a
partition of a total number $n$ of identical particles (animals,
atoms, stars, human beings, etc) into clusters (=groups) of
different sizes, so that $k_j$ is the number of clusters of size
$j.$ The only possible infinitesimal (in time) transitions are
coagulation (merging) of two clusters of sizes $i$ and $j$ into
one cluster of size $i+j$ and fragmentation (splitting) of a
cluster of size $i+j$ into two clusters of sizes $i$ and $j.$
Given a state $\eta\in \Omega_n,$ with $k_i,k_j>0$ for some $1\le
i,j\le n,$ denote by $\eta^{(i,j)}\in\Omega_n$ the state that is
obtained from $\eta$ by the  coagulation of some two clusters  of
sizes $i$ and $j$, and denote by $u_c(\eta,\eta^{(i,j)})$ the rate
of the infinitesimal transition $\eta\rightarrow \eta^{(i,j)}.$
Similarly,
 for
a given state  $\eta\in\Omega_n$ with $k_{i+j}>0,$ let
$\eta_{(i,j)}$ be the state that is obtained from $\eta$ by the
fragmentation of some cluster of size $i+j$ into two clusters of
sizes $i$ and $j,$ and let $u_f(\eta,\eta_{(i,j)})$ be the rate of
the infinitesimal transition $\eta\rightarrow \eta_{(i,j)}.$
Denoting by
$$q(\eta;i,j)=\frac{u_c(\eta,\eta^{(i,j)})}{u_f(\eta^{(i,j)},\eta)}$$
the ratio of the above transitions, the  important property of
reversibility of multiplicative measures is derived by verifying
 the detailed balance condition.

\nin \begin{proposition} A multiplicative measure $\mu_n$ defined
by \refm[3] is reversible with respect to the  transition rates
$u_c,u_f$ if  their ratio satisfies:
 \be
 q(\eta;i,j)=
  \begin{cases}
    \frac{a_{k_i-1}^{(i)}a_{k_j-1}^{(j)}a_{k_{i+j}+1}^{(i+j)}}
{a_{k_i}^{(i)}a_{k_j}^{(j)}a_{k_{i+j}}^{(i+j)}},&\quad
   \text{if}\ \ i\neq j:k_i,k_j>0\\\\
       \frac{a_{k_i-2}^{(i)}a_{k_{2i}+1}^{(2i)}}
{a_{k_i}^{(i)}a_{k_{2i}}^{(2i)}},&\quad \text{if}\ \  i=j: k_i\ge
2.
  \end{cases}
\la{cv}
  \end{equation}

\end{proposition}
\vskip 1cm  An immediate consequence  of  the Proposition 1 is
that a multiplicative measure $\mu_n$ defined by \refm[3] is the
equilibrium distribution of a $CFP$ with transition rates  $u_c,\ u_f$ obeying
the condition \refm[cv], under some sequence of probabilities $\{a_{k_j}^{(j)}\}.$

We now distinguish a class of $CFP's$ with transition rates $u_c,\ u_f$ of the
form  $$ u_c(\eta,\eta^{(i,j)})=
  \begin{cases}
    k_ik_j\phi(i,j) , & \text{if}\quad i\neq j, \quad k_ik_j>0 \\
    k_i(k_{i}-1)\phi(i,i) , & \text{if}\quad i=j, \quad k_i\ge 2,
  \end{cases}
$$
\be u_f(\eta,\eta_{(i,j)}) =k_{i+j}\varphi(i,j), \quad 2\le i+j\le
n,\quad k_{i+j}\ge 2, \la{qw}
 \end{equation}
where $\phi,\varphi $ are some symmetric nonnegative functions on
the set of pairs of positive integers. Treating   the functions $
\phi,\varphi $ in \refm[qw] as the rates of a single coagulation
and a single
 fragmentation respectively, the induced $CFP's$ can be viewed as  mean-field
 models on the set $\Omega_n$. In fact, \refm[qw] tells us that at any state $\eta\in \Omega_n,$
each cluster may coagulate with every other one or may be
fragmented into  two parts, so that the net rates of the
transitions $\eta\rightarrow \eta_{(i,j)}$  and $\eta\rightarrow
\eta^{(i,j)}$ are  sums of the rates of all possible single
coagulations and
 fragmentations respectively, at the state
$\eta$. It is proven in \cite{DGG} with the
 help of the Kolmogorov cycle condition that in the case when all rates of single
  coagulations  and fragmentations
   are positive,
 the  $CFP's$ given by
 \refm[qw] are reversible if and only if the ratios of the single
 transitions are of the form
\be\frac{\phi(i,j)}{\varphi(i,j)}= \frac{a_{i+j}}{a_{i}a_j}, \quad
i,j\ge 1,\la{143}\end{equation} with some  $a_i>0, \quad i\ge 1.$
The corresponding $CFP's$ are known as classical reversible models
of clustering and networks studied in 1970-s by Kelly and Whittle
(see \cite{DGG}-\cite{EG} and references therein). The equilibrium
measures of the mean- field $CFP's$ with the rates \refm[qw],
\refm[143] are multiplicative measures $\mu_n$ induced by the
conditioning relation \refm[1] with $Z_j$ distributed $Po(a_j), \
j\ge 1$. An  example of a reversible  $CFP$ which is not a mean
field, is provided by setting in \refm[2] $Z_j=G_j-1$, where $G_j$ is
distributed geometrically with parameter $p^j,\ 0<p<1$. Then
$a_{k}^{(j)}=p^{jk}q_j, \ q_j=1-p^j, \ j\ge 1, \ k\ge 0 ,$ and the
corresponding measure $\mu_n$ is the  uniform one on the set
$\Omega_n,$ while in \refm[cv], $q(\eta;i,j)\equiv 1.$ It is simple to see  that, due to the last fact, the
 net transition rates of form \refm[qw] do not provide the detailed balance condition for the $CFP$
considered in the example, which implies that the above reversible $CFP$ is not a mean -field model.

 \nin $\bullet$ {\bf $CFP's$ on set partitions}(\cite{BP},\cite{P},
 \cite{Bert}).
We assume here that in the preceding set up for $CFP's$, particles
are labeled by $1,\ldots ,n,$ so that  the state space of the
system of clusters related to a $CFP$, becomes the set
$\Omega_{[n]}=\{\pi_{[n]}\}$ of all partitions $\pi_{[n]}$ of the
set $[n]=\{1,\ldots,n\}$ into subsets. Recall that a partition of $[n]$ into
$k$ blocks (clusters) $A_1,\ldots, A_k$ is
$\pi_{[n],k}=(A_1,\ldots, A_k),$ where $A_j, \ 1\le j\le k\le n$
are nonempty and disjoint subsets of $[n]$ whose union is $[n]$
and which are numbered, e.g. in the order of their least element.
Denoting $\vert A_j\vert$ the size of a cluster $A_j,$ we further
assign to  each $A_j,$ a weight $m_{\vert A_j\vert} $ which is a
number of  possible states of $A_j$, the states can be e.g.,
shapes (in the plane or in space), colors, energy levels, etc.
This says that to the set partition $\pi_{[n],k}$ correspond
$\prod_{j=1}^k m_{\vert A_j\vert} $ different structures with the
same blocks $A_1,\ldots, A_k$, so that the total number of
structures formed by all partitions of the set $[n]$ into $k$
given clusters is equal to \be \sum_{\pi_{[n],k}\in
\Omega_{[n],k}}\prod_{j=1}^k m_{\vert
A_j\vert}:=B_{n,k},\label{Bnk}\end{equation} where $B_{n,k}$ is
known as a Bell polynomial in  weights $m_1,\ldots, m_{n-k+1}.$
Similar to the  setting for decomposable combinatorial
structures,  a random structure $\Pi_{[n],k}$ is the one chosen
randomly from the set of $B_{n,k}$ structures. As a result, for a
given $k$ a measure $p_{[n],k}$ on the set
$\Omega_{[n],k}=\{\pi_{[n],k}\}$ is induced: \be
p_{[n],k}(\pi_{[n]})=\frac{\prod_{j=1}^k m_{\vert
A_j\vert}}{B_{n,k}},\quad
\pi_{[n]}\in\Omega_{[n],k}.\label{gb}\end{equation} In a more
general setting which encompasses a variety of models (see
\cite{BP},\cite{P}), the weights $m_j$ in \refm[gb] are allowed to
be arbitrary nonnegative numbers. Pitman \cite{P} calls the
$\Pi_{[n],k}$ a Gibbs partition and the measure $p_{[n],k}$
microcanonical Gibbs distribution. Obviously, the vector $(\vert
A_1\vert,\ldots \vert A_k\vert)$ of block size counts defines a
partition of the integer $n$ into $k$ summands, induced by the
generic  set partition $\pi_{[n],k},$ and it is known that  to
each $\eta=(k_1,\ldots,k_n)\in \Omega_n,$ such that
$k_1+\ldots+k_n=k,$ correspond
$$\frac{n!}{\prod_{j=1}^n (k_j!)(j!)^{k_j}}$$
different set partitions $\pi_{[n],k}\in \Omega_{[n],k},$ each one of
them having   the same probability $p_{[n],k}(\pi_{[n],k})$ given
by \refm[gb]. Thus, the Gibbs distribution $p_{[n],k}$ on
$\Omega_{[n],k}$ induces the Gibbs distribution $p_{n,k}$ on the
set $\Omega_{n,k}$ of integer partitions of $n$ into  $k$ positive
summands: \be
p_{n,k}(\eta)=(B_{n,k})^{-1}\prod_{j=1}^n(\frac{m_j}{j!})^{k_j}\frac{1}{(k_j)!},
\quad \eta=(k_1,\ldots,k_n)\in \Omega_{n,k}, \label{pnk}
\end{equation}
where the partition function $B_{n,k}$ defined as in \refm[Bnk]
can be rewritten in the following form: \be B_{n,k}=\sum_{\eta\in
\Omega_{n,k}}\prod_{j=1}^n(\frac{m_j}{j!})^{k_j}\frac{1}{(k_j)!}.
\end{equation}
 From \refm[pnk] it is easy to derive   Kolchin's
representation of Gibbs partitions (see \cite{P}, Theorem 1.2). On
the other hand, the distribution $p_{n,k}$ given by \refm[pnk] is
produced by conditioning the multiplicative measure
 $\mu_n$ defined by \refm[3] on the event $Z_1+\ldots
+Z_n=k,$ with $Z_j\sim Po(a_j), \ a_j=\frac{m_j}{j!},\ j\ge 1.$
However, this embedding of the generic model associated with set
partitions of $[n]$ into the setting for conditioning relation
\refm[1] does not facilitate the study of a wealth of problems
(see \cite{BP}) arising from treating $p_{[n],k},\ k=1,\ldots,n$
as marginal distributions of irreversible time continuous markov
processes of pure fragmentation (or pure coagulation) on the state
space $[n]$. The study of these problems  was initiated by Kingman
and Pitman and has been extensively continued by a group of
researchers including Pitman, Bertoin, Berestycki, Gnedin et al.

\vskip .5cm In what follows we will refer to all  models induced
by the conditioning relation \refm[1] as random structures $(RS's)$.

\section{Objective and Summary.}

 \nin In this paper, we  study the asymptotic behaviour, as $n\to\infty,$ of the
 random
 vector
 $(K^{(n)}_1,\ldots,K^{(n)}_l)$ composed of the first $l\ge 1$ components of the random
 vector
${\bf K}^{(n)}$, defined by \refm[1].
  In view of  the independence of the random variables $Z_j, \
j\ge 1$ in  \refm[1], there  was a common belief in physics and
combinatorics that the small (compared with $n$) counts
$K^{(n)}_1,\ldots,K^{(n)}_l$ become independent, as $n\to \infty,$
for any fixed $l\ge 2$, and this was proven in
 a variety of particular cases of $RS's.$ We show that in general
 the assumption of asymptotic  independence fails.
 This said, we note that properly scaled large component
counts $K^{(n)}_l,K^{(n)}_{l+1},\ldots,K^{(n)}_n$ are known to be
dependent in the limit, for any fixed $l\ge 1.$
  Our main result which is Theorem
1 in Section 3, consists of establishing the necessary and sufficient conditions for the asymptotic
independence of small component counts. Combining  this result with  the
Schur's lemma we provide in Section 4   a plausible sufficient condition for
convergence  $RS$'s. This allows us to answer the question of convergence of counting processes for
the three basic types of $RS$'s
 discussed in Section 1. It turns out that many
models of $RS$'s  are divergent. In a parallel way we discuss the
problem of convergence for $CFP$'s. The final section, Section 5,
contains concluding remarks, among them a historical background of
the problem.
\section{Main result.}
\begin{definition}. We say that the counting process $\{{\bf K}^{(n)},\ n\ge 1\},$
 is convergent,  if for each fixed \ $l\ge 1,$ the
probability law $ {\cal L}(K_1^{(n)},\ldots,K_l^{(n)})$ weakly
converges, as $n\to \infty$, to some probability law \ $F_l$  on
$\ql,\ l\ge 1.$ Moreover, we say that the counts
$K_1^{(n)},\ldots,K_l^{(n)}, \ l\ge 2$ of small components of the
random vector\ ${\bf K}^{(n)}$ are
 asymptotically independent  if
 the above laws $F_l$ are product measures on $\ql,$  for all
finite $l\ge 2.$
\end{definition}

\nin Note that in contrast to the setting for limit shapes (see
e.g.\cite{V1}, \cite{EG1}), in this paper we are interested in the
weak convergence of non scaled multiplicative measures.

\nin  For a fixed $l\ge 1$,  given $k_1, \ldots, k_l$ and
sufficiently large $n$, we denote \be M_l=\sum_{j=1}^l jk_j, \quad
and
  \quad T^{(l)}(n, M_l)=\P(\sum_{j=l+1}^{n-M_l} jZ_j=n-M_l). \la{6}
\end{equation}

\nin It is immediate that \be T^{(l)}(n,M_l)=T^{(l)}
(n-M_l,0):=T^{(l)}_{n-M_l}, \quad l\ge 1,\quad M_l\ge 0. \la{7}
\end{equation}

\nin Assuming  in what follows that $$a^{(j)}_0>0, \ j\ge 1,$$ we
will be dealing with the ``scaled" quantities $\tilde a^{(j)}_k$,
$\tilde c_n,$ and $\tilde T^{(l)}_{n-k}$ defined  by \be \tilde
a^{(j)}_k=\frac{ a^{(j)}_k}{ a^{(j)}_0},\quad k\ge 0,\quad j\ge 1,
\label{ak}\end{equation} \be \tilde c_n=\big(\prod_{j=1}^n
a^{(j)}_0\big)^{-1}c_n, \quad n\ge 1, \quad \tilde c_0=1, \la{cnk}
\end{equation}\be \tilde T^{(l)}_{n-k}= \Big(\prod_{j=l+1}^{n-k}
a^{(j)}_0\Big)^{-1}T^{(l)}_{n-k}, \quad \tilde
T^{(l)}_{0}=1,\quad l\ge 1,\quad 0\le k\le n.\la{11}
\end{equation}
In the context of decomposable combinatorial structures,  the
quantities $\tilde{c}_n$ and $\tilde T^{(l)}_{n-k}$ have a
significant combinatorial meaning. Denoting by $p_n$ the number of
structures of size $n$,  we demonstrate in Section 4 that
$$ p_n=
\left\{  \begin{array}{ll}
    p^{-n}\tilde{c}_n, & \hbox{\text{for multisets ($NB(p^j,m_j)$) and selections ($Bi(p^j,m_j)$})} \\
    n!\tilde{c}_n, & \hbox{\text{for assemblies ($Po(a_j)$).}}
  \end{array}
\right.
$$
In analogous way, $\tilde T^{(l)}_{n-k}$ is related to the number
of structures of size $(n-k)$ with all  component sizes greater
than $l$.

  With the help of the above  notation, we have
$$ \P(
K_1^{(n)}=k_1,\ldots,K_l^{(n)}=k_l)=c_n^{-1}\big(\prod_{j=1}^l
a_{k_j}^{(j)}\big)\P(\sum_{j=l+1}^n jZ_j=n-M_l)=
$$
\be \Big(\prod_{j=1}^l \tilde a_{k_j}^{(j)}\Big) \frac{\tilde
T^{(l)}_{n-M_l}}{\tilde{ c}_n}, \la{zyz}
\end{equation}
where in the last step we have used the fact that
$$\P(\sum_{j=l+1}^n
jZ_j=n-M_l)=T^{(l)}_{n-M_l}\prod_{k=1}^{M_l}\P(Z_{n-M_l+k}=0)$$
and the definitions \refm[ak], \refm[cnk] and \refm[11] of the
``scaled" quantities.
 Note that in view of \refm[6],
$\tilde T^{(l)}_{n-M_l}$ is the same for all $k_1,\ldots, k_l:
\sum _{j=1}^l jk_j=M_l.$

 Central to our subsequent study is the notion of smoothly
growing real sequences $RT_\rho,$ the definition of which we adopt from
\cite{Bur}, \cite{BB}.
\begin{definition}.
$RT_{\rho}, \ 0\le\rho\le\infty$ is the collection of sequences
$\{d_n\}_{n\ge 1}$ of \ nonnegative  numbers that satisfy

\be \lim_{n\to\infty}\frac{d_n}{d_{n+1}}=\rho.   \la{9}
\end{equation}
\end{definition} \qed

 \nin   Sequences in $RT_\rho$ play a key role in Compton's theory of logical limit
laws  and in additive number theory (for references see
\cite{Bur}, \cite{B},\cite{BB}).

\nin Now we are prepared to state our main result.

\begin{thm}.
 The counting process $\{{\bf K}^{(n)},\ n\ge 1\}$ is
convergent if and only if the following two conditions hold:

\nin (a) $\{\tilde{ c}_n\}_{n\ge 0}\in RT_\rho,$ for some
$0\le\rho<\infty,$ and

\nin (b) For each $l\ge 1$ there exists a   positive finite limit \be q^{(l)}=\lim_{n\to
\infty}\frac{\tilde{T}_n^{(l)}}{\tilde{c}_n}. \la{el}\end{equation}

Moreover, counts of small components of a convergent random vector ${\bf
K}^{(n)}$ are asymptotically independent.
\end{thm}


\nin {\bf Proof.}

 \nin In view of \refm[zyz], the counting process $\{{\bf K}^{(n)},\ n\ge 1\}$ is convergent  if and only if the fraction on the RHS of \refm[zyz] has finite limits, as $n\to \infty$ for all fixed $M_l\ge 0,\ l\ge 1,$ while for any $l\ge 1$ there exists an $M_l\ge 0,$ such that \be \lim_{n\to
\infty}\frac{\tilde T^{(l)}_{n-M_l}}{\tilde{c}_n}>0.\la{el1}
\end{equation}
We note that \refm[el1] secures that the limiting distribution is a probability measure.
\nin Next we  write  \be\frac{\tilde
T^{(l)}_{n-M_l-1}}{\tilde{c}_n}=\left(\frac{\tilde
T^{(l)}_{n-M_{l}-1}}{\tilde{ c}_{n-1}}\right)\left(\frac{\tilde
c_{n-1}}{\tilde{c}_n}\right).\la{f,}\end{equation}
 We first prove the necessity of the conditions $(a),(b)$ of the theorem. Assuming that  $\{{\bf K}^{(n)},\ n\ge 1\}$
converges, it follows from \refm[zyz] that there exists a finite limit
$$\lim_{n\to \infty}\P(
K_1^{(n)}=0,\ldots,K_l^{(n)}=0)=\lim_{n\to \infty} \frac{\tilde
T^{(l)}_{n}}{\tilde c_n}:=q^{(l)}<\infty, \ l\ge 1.$$
Consequently,\refm[f,], \refm[el1] and the preceding discussion imply  that  $\tilde c_n\in
RT_\rho,$ for some $0\le\rho<\infty$  and  we get from \refm[f,]
\be \lim_{n\to \infty} \frac{\tilde T^{(l)}_{n-M_l}}{\tilde
c_n}=q^{(l)}\rho^{M_l},  \ l\ge 1,\la{bo}\end{equation} for all $
0\le M_l<\infty.$ From the latter and \refm[el1] we conclude that $q^{(l)}$
is positive. For the proof of sufficiency we first apply \refm[f,] with $M_l=0$ to conclude, by virtue  of the conditions (a) and (b), that \refm[bo]
holds with $M_l=1, \l\ge 1$ and so on, proving \refm[bo] for all $M_l\ge 0.$
As a result, \be \lim_{n\to \infty} \P(
K_1^{(n)}=k_1,\ldots,K_l^{(n)}=k_l)=q^{(l)}\prod_{j=1}^l \tilde {a}_{k_j}^{(j)}\rho^{jk_j},\ l\ge 1,\la{ft}\en
  by \refm[zyz] and the definition of $M_l$. Since the sum over $(k_1,\ldots,k_l)\in \ql$ of  the LHS of \refm[ft] is equal to 1, we obtain the explicit expression for $q^{(l)}:$
\be q^{(l)}=\big(\prod_{j=1}^l\tilde{S}^{(j)}(\rho)\big)^{-1},\ l\ge 1, \la{th}\en
where we denoted $\tilde{S}^{(j)}(\rho)=\sum_{k\ge 0}\tilde{ a}_{k_j}^{(j)}\rho^{jk_j},$ in accordance with
\refm[sj]. This shows that the limiting distribution  of the probability law $ {\cal L}(K_1^{(n)},\ldots,K_l^{(n)})$  is the product probability
measure \be \prod_{j=1}^l \frac{\tilde {a}_{k_j}^{(j)}\rho^{jk_j}}{\tilde{S}^{(j)}(\rho)}, \ l\ge 1.\la{mr}\en \qed

\begin{remark}
(i) Setting $0^0=1$ and recalling that $\tilde {a}_{0}^{(j)} =1,\ j\ge 1,$ it follows from \refm[zyz] and \refm[mr]
that in the case of a convergent counting process with $\rho=0$ in condition (a), the limit law of the random probability vector $(K_1^{(n)},\ldots,K_l^{(n)}), \ l\ge 1$ is the measure concentrated at the singleton $(0,\ldots,0)\in \ql.$
We also observe that in this case $q^{(l)}=1,\ l\ge 1,$ in accordance with \refm[ft], while
$\lim_{n\to \infty} \frac{\tilde T^{(l)}_{n-M_l}}{\tilde
c_n}=0, \ l\ge 1,$ for all $M_l>0$.

\nin (ii) Condition (c)  implies $\{\tilde T^{(l)}_{n}\}_{n\ge
0}\in RT_\rho,$ for all \ $ l\ge 1,$ with the same
$0\le\rho<\infty$ as for the sequence $\{\tilde c_n\}_{n\ge 0}.$
This can be seen by writing
$$\frac{\tilde
T^{(l)}_{n }}{\tilde{c}_n}=\left(\frac{\tilde T^{(l)}_{n}}{\tilde
T^{(l)}_{n-1}}\right)\left(\frac{\tilde
T^{(l)}_{n-1}}{\tilde{c}_{n-1}}\right)\left(\frac{\tilde
{c}_{n-1}}{\tilde{c}_n}\right)$$ and then applying the fact that
$0<q^{(l)}<\infty$.
Conversely, if \refm[el] holds and $\{\tilde T^{(l)}_{n}\}_{n\ge
0}\in RT_\rho,$ for some $l\ge 1$ and $0\le\rho<\infty,$ then $\{\tilde c_n\}_{n\ge 0}\in RT_\rho,$
with the same $\rho$. \\
\nin (iii) One can see from the proof of Theorem 1 that the
condition $q^{(l)}>0, \ l\ge 1$ which is a part of the condition
(c), ensures the tightness of the corresponding sequences of
finite dimensional probability measures.
 \end{remark}

To formulate the forthcoming corollary, we need to extend the
definition \refm[ti] of the tilting transformation to the case
$\rho=0,$ in the following natural way: \be a^{(j)}_k(0)=
  \begin{cases}
    1 , & \text{if}\ k=0 \\
    0 , & \text{otherwise}.
  \end{cases} \label{ro}
\end{equation}

\begin{cor} Let ${\bf
K}^{(n)}$ be a convergent counting process, such that  $\tilde
c_n\in RT_\rho,$ for some $ 0\le\rho<\infty.$ Then \be \lim_{n\to
\infty}\P( K_1^{(n)}=k_1,\ldots,K_l^{(n)}=k_l)=\prod_{j=1}^l
a_{k_j}^{(j)}(\rho), \quad l\ge 1, \la{111}
\end{equation}
\nin where $ a_{k_j}^{(j)}(\rho)$  are the generic probabilities \
$ a_{k_j}^{(j)}$ tilted with the above $\rho .$
\end{cor}

\nin {\bf Proof.} By \refm[el]   and  the definitions
\refm[ti],\refm[ro] and  \refm[ak], it follows from Theorem 1 that
for a convergent counting process, $$ \lim_{n\to \infty}\P(
K_1^{(n)}=k_1,\ldots,K_l^{(n)}=k_l)
=\prod_{j=1}^l\frac{\tilde{a}_{k_j}^{(j)}\rho^{jk}}{\tilde{S}^{(j)}(\rho)}=$$

$$
\prod_{j=1}^l\frac{a_{k_j}^{(j)}\rho^{jk}}{S^{(j)}(\rho)}=
\prod_{j=1}^l a_{k_j}^{(j)}(\rho), \ \ l\ge 1.$$
 \qed

\begin{remark}

(i)  Corollary 1 says that asymptotic  independence with the limit
product measure composed of  generic probabilities $
a_{k_j}^{(j)}$ takes place only if $\{\tilde c_n\}\in RT_1.$

(ii)  We denote by $c_n(\theta)$ the quantity $c_n$ corresponding
to the tilting of the probabilities $a_k^{(j)}$ with $\theta\ge 0$
and recall that a multiplicative measure $\mu_n$ is invariant
under all possible tiltings  of the probabilities with $\theta>0.$
By definitions \refm[ti] and \refm[cnk] we then have
$$\widetilde{c_n(\theta)}=\theta^n\tilde{c}_n,\ n\ge 0, $$
where $\widetilde{c_n(\theta)}$ is the scaling of $c_n(\theta).$
So, if $\{\tilde{ c}_n\}_{n\ge 0}, \in RT_\rho,$ for some
$0\le\rho<\infty,$ then
$$\{c_n(\theta)\}_{n\ge 0}\in
RT_{\frac{\rho}{\theta}}.$$ This clarifies the following meaning
of Corollary 1. Consider   a  convergent  counting process such
that $\{\tilde{ c}_n\}_{ n\ge 0}\in RT_\rho,$ with some $\rho>0.$
Then the  whole family of counting processes obtained by tilting
the original one with all possible $\theta>0,$
 has the same limit
finite dimensional distributions as the counting process obtained
by tilting the original one with the  $\rho>0$, so that the
corresponding quantity $\{\widetilde{ c_n(\rho}), \ n\ge 0\}\in
RT_1.$
\end{remark}
\section{ Convergent and divergent random structures.}
We agree to call a $RS$ convergent/divergent if the corresponding
 counting process converges/diverges  in the sense of
Definition 1.

Assuming that condition (a) of Theorem 1 holds, our tool for verifying   condition (b)  for
the
 models considered will be the remarkable Schur tauberian lemma
cited below.
 With an obvious abuse of notation, we say that
a power series $f(x)=\sum_{n\ge 0} d_nx^n $ is in $RT_\rho$ if
$\{d_n\}_{n\ge 0}\in RT_\rho.$ We denote by * the Cauchy product,
which is  extended to formal power series as usual(see
\cite{Bur}).
\begin{lemma}(Schur (1918), see \cite{Bur},p. 62).

Let $f=f_1*f_2,$ where $f, f_1, f_2$ are power series with
 coefficients $d_n, d_n^{(1)}, d_n^{(2)}, \ n\ge 0,$
respectively,

 such that :

 (a) $f_1\in RT_\rho$ for some $0\le \rho<\infty
$ and

(b) the radius of convergence of $f_2$ is greater than $\rho.$

Then \be \lim_{n\to \infty}\frac{ d_n}{d_n^{(1)}}=f_2(\rho).
\la{shu}\end{equation}
\end{lemma}
\vskip.5cm
Schur's lemma is widely used in asymptotic enumeration and in the study of asymptotic densities of additive number systems (see \cite{Bur}). The proof of the lemma is quite simple (see \cite{PS}, Problem 178).
In \cite{Ye} a version of Schur's lemma for Dirichlet series was obtained, which allowed applications to multiplicative number theory.

We  note   the fact that, under the conditions of Schur's lemma, $f_1\in RT_\rho$  implies $f\in RT_\rho.$ This can be seen by writing
$$\frac{d_{n-1}}{d_n}=\Big(\frac{d^{(1)}_{n-1}}{d^{(1)}_n}\Big)
\Big(\frac{d_{n-1}}{d_{n-1}^{(1)}}\Big)
\Big(\frac{d_{n}^{(1)}}{d_n}\Big).$$

We now outline the scheme of application of Schur's lemma to our specific setting.
 Treating $S^{(j)}$ in \refm[sj] as the generating probability function  of the
random variable $jZ_j$ (=of the sequence of probabilities $\{a^{(j)}_k\}$)  in the conditioning
relation \refm[1], we have:
 $$S^{(j)}(x)=\sum_{k\ge 0} a_{k}^{(j)} x^{kj}.$$
\nin  Clearly, the radius of convergence of $S^{(j)}$ is $\ge 1,$
for all $j\ge 1.$ Next, we denote $\tilde{S}^{(j)}=\frac{1}{a^{(j)}_0}S^{(j)}$,
$\tilde{g}=\prod_{j\ge 1}\tilde{S}^{(j)}$ and $\tilde{g}_{\tilde{T}^{(l)}}=\prod_{j\ge l+1}\tilde{S}^{(j)}.$
Since $\tilde{a}_{0}^{(j)}=1,\ j\ge 1,$
the above are the generating functions for the scaled sequences $\{\tilde{a}_{k}^{(j)}\}_ {k\ge 0},\
\{\tilde{c}_n\}_{n\ge 0}$ and $\{\tilde{T}^{(l)}_n\}_{n\ge 0}$
respectively, as defined by \refm[ak]-\refm[11].
 Finally, writing  $\tilde g^{(l)}=\prod_{j=1}^l \tilde
S^{(j)},\ l\ge 1$ we have

\be \tilde g=\tilde g_{\tilde T^{(l)}}*\tilde
g^{(l)}\la{112a}\end{equation}

\nin or, equivalently,
 \be \tilde g_{\tilde T^{(l)}}=\tilde
g*\left(\frac{1}{\tilde g^{(l)}}\right).\la{112}\end{equation}

\begin{remark}
We make here use of the representation \refm[112a] to show that condition (a) of Theorem 1 does not imply even the  existence of the limit \refm[el] defining  $q^{(l)}, \ l\ge 1.$
Let $l=1$, $$\tilde{g}(x)=\frac{1}{1-x}, \quad \tilde{g}^{(1)}(x)=1+x,\quad  \tilde{g}_{\tilde{T}^{(l)}}(x)=\frac{1}{1-x^2},\quad \vert x\vert<1.$$ Then $$\tilde{c}_n=1, \quad n\ge 1,
\quad \tilde{T}_n^{(1)}=\left\{
                          \begin{array}{ll}
                            0, & \hbox{if}\ \ n \ \text{is \ odd}\\
                            1, & \hbox{if}\ \ n\ \text{is \ even}
                          \end{array}
                        \right. \ n\ge 1,
$$
which shows that the $\lim_{n\to \infty}\frac{\tilde{T}_n^{(1)}}{\tilde {c}_n}$ does not exist.
The scheme considered is realized by the following sequence of random variables $Z_j,\ j\ge 1:$
$$Z_1\sim Bernoulli(1/2),\quad Z_j\sim Po(a_j), \ \text{with}\ \ a_j=\left\{
                                                        \begin{array}{ll}
                                                          0, & \ \hbox{if}\ \ j\ \text{is \ odd} \\
                                                          2/j, & \hbox{if}\ \ j\ \text{is \ even}
                                                        \end{array}
                                                      \right..
$$
\end{remark}

\vskip .5cm
 \begin{proposition}( {\bf Sufficient condition of convergence/divergence}).

 Let $\tilde{g}\in RT_{\rho},$ $0\le \rho<\infty$ and let the radius of convergence of the series  $\frac{1}{\tilde{g}^{(l)}}$
 be greater than $\rho$, for all $l\ge 1.$ Then a $RS$ converges if \ $\Big(\tilde{g}^{(l)}(\rho)\Big)^{-1}>0,\ l\ge 1$ and it diverges if \
$\Big(\tilde{g}^{(l)}(\rho)\Big)^{-1}=0, \ l\ge 1.$
\end{proposition}
{\bf Proof}. Applying Schur's lemma to \refm[112], we get $q^{(l)}=\Big(\tilde{g}^{(l)}(\rho)\Big)^{-1}<\infty.$
By Theorem 1 this implies the claim.\qed

\begin{remark}
The example in Remark 3 demonstrates the importance of the second condition of the Proposition 2. In fact, in the model considered $\tilde{g}\in RT_{1}$ and $\big(\tilde{g}^{(1)}(1)\big)^{-1}>0.$ However the $RS$ diverges, since the radius of convergence of
\ $\big(\tilde{g}^{(1)}(x)\big)^{-1}=\frac{1}{1+x}$ equals to $1$.

\end{remark}

Proposition 2 allows to suggest the following  two-step strategy for deciding about\\ convergence/divergence of $RS$'s.

 (i) Validation of
the condition $\tilde g\in RT_\rho$ for some $0\le\rho<\infty.$
Our treatment of the problem  is based on application of known
 sufficient conditions on  sequences $\{m_j\}_{  j\ge 1}$ that
guarantee the
  $RT_\rho$ property for  the induced sequences $\{\tilde{ c}_n\}_{ \ n\ge 0}. $ The conditions  we employ are the products  of two quite
different lines of research: \begin{itemize} \item Burris-Bell theory
(\cite{B}, \cite{BB}) of $RT_\rho $ sequences, developed  with the
help of analytical tools stemming from Tauberian theory.
Motivation of this research came from
 Compton's (1980s) theory of logical limit laws and also from the additive
number system theory.

\item
Sufficient conditions for $\tilde g\in RT_\rho$ implied by asymptotic formulae for the number of decomposable structures, in
particular recent results  by Barbour, Freiman, Stark and the
author, derived with the help of probabilistic methods.
\end{itemize}
 Each one
of the two approaches has its particular limitations, some of
which are noted later on in this Section. In this connection we
observe that there is little hope to obtain plausible necessary
conditions on
$\{m_j\}_{ j\ge 1}$ implied by the $ RT_\rho$ property  of \ $\tilde g$. \\

(ii) Validation of the condition (b) of Schur's lemma for the functions $(\tilde{g}^{(l)})^{-1},\ l\ge 1.$
Provided  the condition holds, the conclusion regarding convergence/divergence is based on the claim of \\
Proposition 2.
\vskip.5cm

Following the aforementioned  strategy, we  examine now the convergence of counting
processes for the three basic types of $RS$'s: assemblies, multisets
and selections, described  in Section 1. Furthermore, the results
obtained explain the crucial difference in the asymptotic
behaviour of mean-field and non mean-field $CFP$'s associated with
the above structures.

First, following  \cite{bar} we will say that a $RS$ is   regularly
varying in case it is
  induced by random variables $Z_j, \ j\ge 1$ in
\refm[1], with \be \E Z_j\sim const \ y^j j^{\alpha}, \quad \alpha
\in \R,\quad y>0,\quad j\to \infty.\la{var}\end{equation} Since
the asymptotic behaviour of regularly varying $RS$'s appears to be
in accordance with the behaviour of the series $\sum_{j=1}^\infty
j^\alpha$ (see \cite{ABT}, \cite{frgr1}, \cite{bar}, \cite{GS}),
 it was suggested in \cite{bar} to
distinguish  the following three classes of regularly varying
structures: logarithmic ($\alpha=-1$), convergent
  ($\alpha<-1$) and
expansive ($\alpha>-1$). (It goes without saying that in this
classification, the meaning of a convergent structure is different
from the one in the present paper). As in \cite{frgr1}, we extend
the above definition of expansive structures to include random
structures with $\E Z_j, \ j\ge 1$ oscillating (in $j$) between
two regularly varying functions,  namely
$$(\E Z_j, \ j\ge 1)\in \calF(r_1,r_2;y):=\{f=f(j):  \gamma_1y^jj^{r_1-1}\le f(j) \le
\gamma_2
 y^jj^{r_2-1}, \ j\ge 1, \ y>0 \},
$$
\nin where  $\gamma_i,\ i=1,2$ are positive constants and $0<
r_1\le r_2.$ (The  requirement $r_1,r_2>0$ is  the characteristic
feature of the expansive case).

 \nin {\bf  Assemblies} \ \ Let
$Z_j\sim Po(a_j)$, \ \  $a_j>0,\ \ j\ge 1$ . In this case,

$$\tilde a_{k_j}^{(j)}=\frac{a_j^{k_j}}{k_j!},\ \
\tilde S^{(j)}(x)=\exp{(a_jx^j)},\ \ j\ge 1$$ and $$\tilde g(x)=
\exp\left(\sum_{j\ge 1}{a_jx^j}\right),\quad \tilde{ g}^{(l)}(x)=
\exp\left(\sum_{j= 1}^l{a_jx^j}\right),\quad l\ge 1,\quad \E Z_j=a_j, \quad j\ge 1.
$$

\nin Thus, the radius of convergence of $\frac{1}{\tilde g^{(l)}}
$ equals  $\infty$ for all finite $l\ge 1.$ Consequently, assuming $\tilde g\in
RT_\rho,\ 0\le\rho<\infty,$ we have $\Big(\tilde{g}^{(l)}(\rho)\Big)^{-1}(\rho)>0.$ By virtue of
Theorem 1, the latter leads to the following

\nin \begin{proposition}  An assembly converges if and only if the
sequence $\{a_j\}_{j\ge 1}$ is such that $\tilde g\in RT_\rho,$
for some $0\le \rho<\infty.$
\end{proposition}
\nin \begin{cor}  Assemblies with the following parameter
functions $a=\{a_j,\ j\ge 1\}$ are convergent:

(i) Smoothly growing: $a\in RT_\rho, \ 0<\rho<\infty$;

(ii) Oscillating: $a\in \calF(\frac{2r}{3}+\epsilon,r;y),$
 where $r,y>0$, while $0<\epsilon\le \frac{r}{3}.$
\end{cor}
\nin {\bf Proof} (i) follows from the important Corollary 4.3 in
\cite{BB}, which says that $a\in RT_\rho, \ 0<\rho<\infty$ implies
$\tilde g\in RT_\rho,$ with the same $\rho.$ For the proof of
(ii), we derive from the asymptotic formula (4.99) in \cite{frgr1}
that for the oscillating assembly
 with $y=1,$ we have $\tilde g\in RT_1.$ Hence, it is left to apply
  (ii) of Remark 2 with $\theta=y$.  \qed

 Note that the
restriction on $\epsilon$ in the part (ii) of the last corollary
determines a bound on the ``size" of oscillation of the function
$\E Z_j=a_j, \ j\ge 1$ that ensures the $RT_\rho$ property for
$\tilde g$. \vskip .5cm
 \nin \emph{Examples.} In combinatorics (see Table 2.2 in
\cite{ABT}),  many assemblies, e.g. permutations
($a_j=\frac{1}{j}$), Ewens sampling formula
($a_j=\frac{\theta}{j}, \ \theta >0$), forests of labeled rooted
trees ($a_j=\frac{j^{j-1}}{j!}\sim const\  e^jj^{-\frac{3}{2}}$),
etc. are regularly varying, which says that in all these cases
$a\in RT_{\frac{1}{y}},\ y>0,$ where $y$ is as in the definition
\refm[var] of a regularly varying structure. In statistical
mechanics, regularly varying assemblies with $y=1$ and $\alpha>0$
are called generalized Maxwell-Boltzman statistics (\cite{V1}).
 By virtue of the condition (i) of Corollary 2, all the
 aforementioned assemblies converge.

Based on  Proposition 3, we give now two examples of divergent
assemblies. Firstly, set partitions ($a_j=(j!)^{-1}, \ j\ge 1$)
diverges, since in this case $\tilde g(x)=e^{e^x-1},$ so that the
radius of convergence of $\tilde g$ is infinity. In this
connection note that the following sufficient condition for the
$RT_\infty$ property of $\tilde g$ was recently established in
\cite{BK}:  If the parameter function $a$ is such that
$gcd\{j:a_j>0\}=1$ and $a_j=O(j^{\theta j}/j!), \ 0<\theta<1,\
j\to \infty,$ then $\tilde g\in RT_\infty$, which means
that the induced assembly diverges.\\
For our second example we construct a divergent assembly with
$\tilde g$ that does not belong to any class $RT_\rho, \ 0\le
\rho\le \infty.$ It is clear that the corresponding sequence
$\{a_j\}_{j\ge 1 }$ should exhibit a wild behavior. We set
$$ a_j=
  \begin{cases}
    j^{-1}, & \text{if $j\ge 1$ is odd}, \\
   j^{-1}+2^{j+1}j^{-1},  & \text{if $j\ge 2$ is even}.
  \end{cases}
$$
\nin We then have
$$\tilde g(x)=\Big(\frac{1}{1-x}\Big)*\Big(\frac{1}{1-4x^2}\Big),$$
\nin  which by the Cauchy product formula gives

$$\tilde c_{n}=\sum_{k=0}^{[n/2]} 4^k= \frac{4^{1+[n/2]}-1}{3},\ \ n\ge 1.$$

\nin Consequently,
$$ \lim_{n\to \infty}\frac{\tilde c_{2n-1}}{\tilde c_{2n}}=1/4$$
$$ \lim_{n\to \infty}\frac{\tilde c_{2n}}{\tilde c_{2n+1}}=1.$$

\nin Finally, note that  for graphs on $n$ vertices,
 $\tilde {c}_n=\frac{2^{n\choose 2}}{n!}$, so that $\tilde{c}_n\in RT_0$ and  by
 (i) of Remark 1, for any $l\ge 1$
 the limit measure is  concentrated on the
singleton $(0,\ldots,0)\in \R_l, \ l\ge 1.$

\vskip.5cm We shift  now to considering reversible   $CFP's$ related
to assemblies. In this case the ratio of the net transitions
\refm[cv] has  the following form:

\be
 q(\eta;i,j)=V(k_i,k_j,k_{i+j})
  \begin{cases}
    \frac{a_{i+j}}{a_{i}a_j}, & \text{if}\ \ i\neq j:k_ik_j>0 \\
    \frac{a_{2i}}{a_i^2} , & \text{if}\ \ i=j:k_i\ge 2,
  \end{cases}\la{ha}
\end{equation}
\nin where \be V(k_i,k_j,k_{i+j})=\begin{cases}
\frac{k_ik_j}{k_{i+j}+1},& \text{if}\ \ i\neq j:k_ik_j>0 \\
\frac{k_i(k_i-1)}{k_{2i}+1},& \text{if}\ \ i=j:k_i\ge 2.
\end{cases}
\end{equation}

\nin As we explained in Section 1, such ratios correspond to
reversible  mean-field $CFP's$ with net transition rates of
coagulation and fragmentation as given by \refm[qw],\refm[143] and
the equilibrium measure $\mu_n$ defined by \refm[3] with $Z_j\sim Po(a_j),\ j\ge 1$. Consequently,
the preceding discussion reveals that amongst mean-field $CFP's$
both convergent and divergent models exist, 
depending on the sequence of parameters $\{a_j\}_1^\infty.$ \vskip.5cm

 \nin
{\bf Multisets} \ Assuming  $Z_j\sim NB(m_j,p^j),$ $0<p<1,$
$m_j\ge 1,$ $j\ge 1,$ we have
$$a_k^{(j)}=(1-p^j)^{m_j}\binom{m_j+k-1}{k}p^{jk},\ k\ge 0,
$$ where $p$ is a free parameter.
This gives $$ \tilde{a}_k^{(j)}=\binom{m_j+k-1}{k}p^{jk},\quad
\tilde{S}^{(j)}(x)=(1-(px)^j)^{-m_j}, \ \vert xp\vert <1,$$ which
leads to the Euler type generating function

\be\tilde{g}(x)=\prod_{j\ge
1}\left(1-(px)^j\right)^{-m_j}, \quad \vert x\vert<p^{-1}.\la{eu}\end{equation}
Clearly, the radius of convergence of $\tilde
g(x)$ is no greater than $p^{-1}.$ Moreover, it is known
(see e.g. \cite{Bur}, Lemma 1.15) that $\tilde g$ converges at
some point $x:\vert xp\vert <1$ if and only if
the sequence $\{m_j\}_1^\infty$ is such that $\lim_{j\to \infty} \left(1-(px)^j\right)^{-m_j}=1.$  Next, in
\refm[112] the function $$\frac{1}{\tilde
g^{(l)}(x)}=\prod_{j=1}^{ l}\left(1-(px)^j\right)^{m_j}>0, \quad \text{for all} \ \ \vert x\vert<p^{-1},\quad l\ge 1.$$

\nin \begin{proposition} A  multiset is convergent if and only if
$\tilde g\in RT_\rho,$ for some $\rho<p^{-1}.$
\end{proposition}
\nin {\bf Proof.} In the case $\tilde g\in RT_{p^{-1}}$,
 we have $ \frac{1}{\tilde
g^{(l)}(p^{-1})}=0,\ l\ge 1,$ which implies divergence. If now $\tilde g\in
RT_\rho$ with $0\le \rho<p^{-1},$ then a multiset converges, by Proposition 2.\qed

\begin{cor} Multisets with the following parameter functions
$m=\{m_j, \ j\ge 1\}$ are divergent: \\ \nin (i) $m\in RT_1;$ (ii)
$m_j=O(j^\alpha),$ for some $\alpha\in R;$
 (iii) $m\in \calF(\frac{2r}{3}+\epsilon,r;1),\ 0<\epsilon\le \frac{r}{3}, \ r>0,$
The aforementioned convergence /divergence hold under any $0<p<1.$

\nin whereas multisets with  parameter functions (iv)- (vi) below
converge:

\nin (iv) $m\in RT_\rho,$ for some $0<\rho <1;$\ \ \
 (v) $m\in \calF(\frac{2r}{3}+\epsilon,r;y), \ 0<\epsilon\le \frac{r}{3},
  \ y>1,\ r>0;$\\
 (vi) $m_j\asymp y^jj^{-1}, y>1. $
\end{cor}
 \nin {\bf Proof.} Each one of the conditions (i)-(v) is
sufficient
 for $\tilde g\in RT_\rho$ with a corresponding $\rho.$ Namely, the conditions
(i) and (iv) are due to  Bell-Burris Theorems 6.1 and 6.2 of
(\cite{BB}), which state that if $m\in RT_\rho,$ for some
$0<\rho\le 1,$ then $\tilde g(p^{-1}x)\in RT_\rho$ with the same
$\rho,$ so that  $\tilde{g}(x)\in RT_{\rho p^{-1}}.$
(Note that
 in Theorem 6.1  in  \cite{BB},  condition (c) is  required only for
the second part of the claim).
 Condition (ii) provides $\tilde g(p^{-1}x)\in RT_1,$ by the
powerful result of Bell \cite{B} that generalizes the Bateman and
Erd$\ddot{o}$s theorem. So, under this condition, $\tilde{g}(x)\in
RT_{ p^{-1}}.$    Conditions (iii) and (v) result from Corollary 2
of \cite{GS} for expansive multisets, which says that in  both
cases $\tilde{g}(p^{-1}x)\in RT_{y^{-1}},$ which is equivalent to
$\tilde{g}(x)\in RT_{y^{-1}p^{-1}}$, with $y\ge 1.$ Regarding the
condition (vi), we firstly recall that $a_j\asymp b_j$ means that
the ratio ${a_j\over b_j},\ j\ge 1$ is bounded above and below by
positive constants.   The sufficiency of the condition (vi) for
$\tilde g\in RT_{y^{-1}p^{-1}}, \ y>1$ was proven by Stark in
\cite{ST},which is  devoted to logical limit laws for logarithmic
multisets. \qed

\vskip.5cm
 \emph{Examples.} Integer partitions
($m_j=1, \ j\ge 1),$  planar partitions ($ m_j=j, \ j\ge 1$) (see
\cite{A}) and generalized Bose-Einstein statistics ($m_j=j^\alpha,
\ \alpha>0$) diverge, since in all these cases $m\in RT_1.$ The
following logarithmic multisets (see \cite{ABT}) with $m\in
RT_\rho, \ \rho<1$ converge : mapping patterns
($m_j\sim \frac{\rho^{-j}}{2j}, \ \rho<1,\ j\ge 1),$ monic polynomials over
$GF(q), \ q>1$ ($m_j\sim \frac{q^j}{j}).$ Also,  forests of
unlabeled, unrooted ($m_j\sim const \ \rho^{-j}j^{-5/2}, \ \rho<1)$ and
rooted ($m_j\sim const \
 b^{-j}j^{-3/2},\ b<1)$ trees converge.
\vskip.5cm \begin{remark}  There is a formal  linkage between
assemblies and multisets,  expressed as follows. In the case of
assemblies the generating function $\tilde g$ is of the
exponential form: $\tilde g(x)=e^{Q(x)},$ where $Q(x)=\sum_{j\ge
1}a_j x^j,$ with $a_j\ge 0,\ j\ge 1.$  For multisets, the Euler
type generating function $\tilde g$ in \refm[eu], can  be written
in the same form, with the function $Q^*$called the star
transformation (see \cite{B}, \cite{BB}) of the generic generating
function $Q(x)=\sum_{j\ge 1}m_jp^jx^j$ of the sequence $\{m_j\}_{
j\ge 1}:$ \be Q^*(x)=\sum_{j\ge 1}m^*_jx^j,\quad
m^*_j=\sum_{lk=j}\frac{m_l p^l}{k}, \quad j\ge 1, \quad
m^*_0=0.\la{star}\end{equation} Thus, $m^*_j\ge 0, \ j\ge 1,$
which is a basic assumption in the theory of $RT_\rho$ sequences.
It was proven in \cite{BB} that if $\{m_jp^j\}_{ j\ge 1}\in
RT_\rho,$ with some $0<\rho<1,$ then $m_j^*\sim m_jp^j,\
 j\to \infty,$ which means that in this case  the  multiset
 behaves asymptotically as the assembly induced by
$Z_j\sim Po(m_jp^j), \ j\ge 1.$ This fact explains the nature of
the condition
 (iv) in Corollary 3.

\end{remark}

Regarding the $CFP$'s associated with multisets,
 \refm[cv] becomes:

\be q(\eta;i,j)=V(k_i,k_j,k_{i+j})
  \begin{cases}
   \frac{m_{i+j}+k_{i+j}}{(m_i+k_i-1)(m_j+k_j-1)} , & \text{if}\ \ i\neq j:k_ik_j>0 \\
  \frac{m_{2i}+k_{2i}}{(m_i+k_i-1)(m_i+k_i-2)} , & \text{if}\ \ i=j:k_i\ge
  2,
  \end{cases}\la{hav}
\end{equation}

\nin where $V(k_i,k_j,k_{i+j})$ is as in \refm[ha].
   The second factor in \refm[hav] depends both on
$\eta=(k_1,\ldots,k_n)$ and the parameters $m_j,\ j\ge 1$ of the
$CFP$, so that the representation \refm[qw] does not hold, which says that  the process is not a mean-field model. To illustrate this
fact, recall that  in a particular case of the $BE$ model ($m_j=1,\
j\ge 1$) we saw in Section 1 that the corresponding  $\mu_n$  is
the uniform measure on $\Omega_n,$ so that we have from \refm[hav]
$$ q(\eta;i,j)=1,\ \eta\in \Omega_n,\ i,j\ge 1,\ i+j\le n. $$

\nin {\bf Selections.} In this case  $Z_j\sim
Bi(\frac{p^j}{1+p^j};m_j),\ m_j\ge O(1),\ j\to \infty$ and $0<p\le 1.$
 Hence, $$\tilde{a}_k^{(j)}=\binom{m_j}{k}p^{jk},$$$$\tilde g(x)=\prod_{j\ge
 1}\left(1+(px)^j\right)^{m_j}.$$ So, similar to the  the case of multisets,
  $\tilde g$ converges at
some point $x:\vert xp\vert <1$ if and only if
the sequence $\{m_j\}_1^\infty$ is such that $\lim_{j\to \infty} \left(1+(px)^j\right)^{m_j}=1,$ which is equivalent
to \be \lim_{j\to \infty}(px)^jm_j=0.\la{dr}\en
Clearly, the radius of convergence, say $\rho$, of $\tilde{g}$ is $\le p^{-1},$ for any nonnegative sequence $\{m_j\}_{j\ge 1}.$ Moreover, the condition   implies that $\rho>0$ for all $0<p\le 1$ if and only if
a selection is expansive (see\cite{bar},\cite{GS}), i.e. $m_j=O(j^\alpha), \alpha>0.$
A majority of practical selections are expansive.

\nin \begin{proposition}  Selections with radius of convergence $0\le\rho<p^{-1}$ converge if and only if
the sequence $\{m_j\}_{j\ge 1}$ is such that $\tilde g\in
RT_\rho.$ All expansive selections converge.
\end{proposition}
\nin {\bf Proof:}
In the case of selections, the radius of convergence of the functions  $$\Big(\tilde{g}^{(l)}(x)\Big)^{-1}=\prod_{j= 1}^l{\big(1+(px)^j\big)^{-m_j}},\ l\ge 1$$ is equal to $p^{-1},$ for all $m_j\ge 0,\ j=1,\ldots,l$, while  $\frac{1}{\tilde{g}^{(l)}(x)}>0, \vert x\vert<p^{-1}.$   Consequently, if  $\tilde{g}\in RT_\rho$ with $\rho<p^{-1},$ then convergence holds by Proposition 2.
However, if $\rho=p^{-1}$ the second condition of Proposition 2 fails, which requires to employ an argument  specific for expansive selections. In this latter case $\tilde{c}_n\to \infty$ according to the asymptotic formula derived in
\cite{GS}. The formula also says that for expansive selections, $\tilde{g}\in RT_{p^{-1}}.$ Moreover, following the   proof of the formula in \cite{GS} it is easy to see that also $\tilde{g}_{\tilde{T}^{(l)}}\in RT_{p^{-1}}.$
   Writing $\tilde{g}(x)=(1+px)^{m_1}\tilde{g}_{\tilde{T}^{(1)}},$ we get $$\tilde{c}_n=\sum_{k\ge 0}\binom{m_1}{k}p^k \tilde{T}^{(1)}_{n-k}, \quad n=1,2,\ldots.$$  With the help of the aforementioned properties of the sequences $\{\tilde{c}_n\}$ and
  $\{\tilde{T}^{(1)}_{n}\}$ we are now  able to write:
  $$\lim_{n\to \infty}\frac{\tilde{c}_n}{\tilde{T}^{(1)}_{n}}=2^{m_1}=
  \tilde{g}_{\tilde{T}^{(l)}}(p^{-1})=\Big(q^{(1)}\Big)^{-1}.$$
  In an analogous way we get $q^{(l)}=\prod_{j=1}^l2^{-m_j}=\Big(\tilde{g}_{\tilde{T}^{(l)}}(p^{-1})\Big)^{-1},\quad l\ge 1.$
 \qed

 For selections,  conditions of
Bell-Burris type on $\{m_j\}_{j\ge 1},$  providing $\tilde g\in
RT_\rho$ are not known. This is due to  the fact that for
selections  the star transformation, as defined by \refm[star],
does not ensure the nonnegativity of $m_j^*, \ j\ge 1.$
However, the  probabilistic method for enumeration of decomposable
structures works in the case of selections also. It follows from
Theorem 5 in \cite{GS}, obtained by implementing the method that
for expansive selections, $m\in \calF(\frac{2r}{3}+\epsilon,r;y),
\ r,y>0,$ implies that $\tilde g\in RT_\rho,$ with
$\rho=(yp)^{-1}.$ By  our Proposition 4, this provides convergence
of the corresponding selections, if $y\ge 1$, under any $0<p<1$.
As a result, we derive the convergence of the following expansive
selections: integer partitions into distinct parts ($m_j=1,\ j\ge
1$) and generalized Fermi statistics ($m_j=j^\alpha,\ \alpha>0).$
In this connection it is in order to note that multisets and
selections with $m_j\equiv 1$ induce uniform measures $\mu_n$ on
the set of integer partitions of $n$ and on the set of integer
partitions of $n$ into distinct parts, respectively. In the first
case the random structure diverges, whereas in the second case
convergence to a Bernoulli product measure holds.

For the associated CFP's we obtain from \refm[cv] \be
q(\eta;i,j)=V(k_i,k_j,k_{i+j})
  \begin{cases}
   \frac{m_{i+j}-k_{i+j}+1}{(m_i-k_i-1)(m_j-k_j-1)} , & \text{if}\ \ i\neq j:k_ik_j>0,\
   0\le k_i\le m_i \\
  \frac{m_{2i}-k_{2i}}{(m_i-k_i+1)(m_i-k_i+2)} , & \text{if}\ \ i=j:2\le
  k_i\le m_i.
  \end{cases}\la{hav1}
\end{equation}

 This shows that, as in the case  of multisets,  the above  CFP's are not mean field models.

\section{ Concluding remarks and history}

Generally speaking, the phenomenon of asymptotic independence  of
a finite number of small groups of particles  in large random
systems (i.e. systems formed of  a large number of randomly
interacting particles) was observed in different fields of
applications, under various mathematical settings. The assumption
of asymptotic independence, sometimes accepted without proof, was
of great help for the study of the probabilistic models
considered. Not pretending to provide a comprehensive survey of
the subject, we point out below a few settings parallel (in some
sense) to the one in the present paper.

\nin (i){\bf The Gibbs conditioning principle} (see \cite{DF},
\cite{DZ}, \cite{GZ}). In the context of an ideal gas model,  the
simplest version of the principle in the title reads as follows. Let $X_1,
X_2,\ldots$ be independent and identically distributed random
variables viewed as energies of individual particles, so that
$X_1+\ldots +X_n$ is the total energy of a system of $n$
particles. Let $\E X_1=1$ and assume some suitable regularity
conditions on a common probability law $P$ of the sequence of
random variables. Then, for a fixed $k\ge 1$ the distribution law
$\Lambda_{k,n}$ of $(X_1,\ldots, X_k)$ conditioned on $X_1+\ldots
X_n=n$ converges weakly, as $n\to \infty$ to the $k$-fold product
law $P^k.$ In statistical physics the law $\Lambda_{k,n}$ is
called a microcanonical distribution, and the gibbs conditional
principle asserts asymptotic independence of energies  of any
finite number of particles in microcanonical ensembles. Formulated
in the beginning of the $20$-th century, the principle has been
extended and refined in different directions, with particular
attention
being paid to the rate of convergence to  limit distributions. \\
To distinguish from the conditioning relation \refm[1],  the
measure $\Lambda_{k,n}$ is defined on the simplex. This fact
implies  that the Cauchy product relationship \refm[112a] for
generating functions, which is basic for the study of
multiplicative measures $\mu_n$, is not valid in this case.
However, we believe that the interplay between the above two
settings deserves further study. Quite independently, the
distribution $\Lambda_{k,n},$ with $k=n$ and discrete and not
necessary identically distributed random variables $X_1,X_2,
\ldots$ was introduced by Kolchin (\cite{Kol}) to  represent the
distributions of cell counts in combinatorial urn schemes. In
\cite{Kol} the representation is called the generalized scheme of
allocation, whereas in \cite{BP} and \cite{P}, it is named the
Kolchin representation formula. We note that the problem of
asymptotic independence is not addressed in \cite{Kol}.

 \nin (ii) {\bf Random combinatorial structures}.
In the theory of random structures, the asymptotic independence of
counts of small components was discussed in numerous papers,
starting from the 1940's. A general set up leading to asymptotic
estimation of the total variation distance between component
spectra of small counts (as defined by the conditioning relation
\refm[1]) and the independent process was developed by Arratia and
Tavar$\acute{e}$ in their seminal paper \cite{AT} (see also
\cite{ABT},\cite{fri} and \cite{Kol}). As a result, asymptotic
independence was established for  logarithmic random structures
with $y=1$ in \refm[var].
In \cite{bar}, the same was proven for regularly varying
convergent ($\alpha<-1$) structures and in \cite{frgr1} the asymptotic
independence was proven  for expansive assemblies with any $y>0$.
Regarding assemblies, multisets and expansive selections,  the
aforementioned results easily follow from our results in Section
4. In fact, recall that for assemblies $\E Z_j=a_j,\ j\ge 1,$  for
multisets $\E Z_j=\frac{m_jp^j}{1-p^j}, \ j\ge 1,\ 0<p<1$ and for
selections $\E Z_j=\frac{m_jp^j}{1+p^j}, \ j\ge 1, \ 0<p<1.$ By
the definition \refm[var] of a regularly varying $RS$ with
parameters $y>0,\alpha\in \R$ we thereby conclude that the
following facts hold.
\begin{itemize} \item For assemblies:  $a=\{a_j, \ j\ge 1\}\in
RT_{y^{-1}}, \ y>0,$ for all $\alpha\in \R,$ which implies
convergence by  condition (i) of Corollary 2.

\item For multisets: $m=\{m_j,  \ j\ge 1\}\in RT_{y^{-1}p}, \
0<p<1,$ for all $\alpha\in R.$ This implies   convergence if $y\ge
1,$ by the condition (iv) of Corollary 3. Note that in the case
$y=p,$ a regularly varying multiset diverges for all $\alpha\in
\R,$ by  condition (i) of Corollary 3.

\item For selections: $m_j\sim const j^\alpha (yp^{-1})^j, \ j\ge
1.$ By the discussion following Proposition 4, this provides
convergence for all $\alpha>0$ and     $y\ge 1 1.$

\end{itemize}

 In this connection, we mention
 that to our knowledge no examples of $RS$'s for which the independence principles fails were given
in the literature, prior to this paper.

 Our results reveal also a basic difference
between  pictures of asymptotic clustering of components in
 convergent $RS$'s and in regularly varying divergent
multisets. Namely, in the case of a convergent $RS$ with $\rho\neq
0$, Corollary 1 tells us that with a positive limit probability
there are components of any fixed sizes, i.e.
$$\lim_{n\to\infty}\P(K^{(n)}_1=k_1,\ldots,K^{(n)}_l=k_l)>0,\quad
k_1,\ldots,k_l\ge 0,\quad l\ge 1.$$ On the other hand, if a
regularly varying multiset diverges, then  $  q^{(l)}=0, \ l\ge
1,$ by Proposition 3, from which it follows that
$\lim_{n\to\infty}\P(K^{(n)}_j=k_j)=0,$ for all finite $k_j\ge 0,
\ j\ge 1.$

\vskip.5cm (iii) {\bf $CFP$'s}. It is common to trace the beginning
of rigorous mathematical models of coagulation-fragmentation back
to the paper by Smoluchowski (1918) where the famous system of
coagulation equations describing the time evolution of the process
was derived. Already in this paper the assumption of independence
(more precisely, the absence of correlations)  of clusters of
small sizes was adopted. Subsequently,  deterministic and
stochastic versions of the model were studied in numerous papers
in probability and various applied fields. The study of reversible
$CFP$'s was concentrated on what we call in the present paper
mean-field $CFP$'s (see \cite{W}, \cite{K}). (Recall that these models
conform to assemblies). In \cite{DGG} the model was treated as a
reversible Markov chain on the set of partitions and it was proven
(Theorem 4, (4.24)) that, if  $\tilde g\in RT_\rho, \
0\le\rho<\infty,$ then at the equilibrium of the process,
$$cov(K_l^{(n)},K_m^{(n)})\to 0,  \ n\to \infty,$$ for any fixed
$l\neq m.$ This is, of course, a weak form of our Proposition 3.
More details on the history of $CFP$'s can be found in \cite{frgr1}.

(iv)  {\bf Convergence of scaled counting processes}.\\
In  the setting of the present paper, as well as in the all above
mentioned settings, the convergence of generic (=nonscaled)
counting processes was studied. Our Theorem 1 asserts that a
generic spectrum either converges (in distribution) to a
distribution with independent components, or diverges. If the
first option is the case, then the simple discrete limiting
process provides approximation of the discrete generic one. In
some (but not all!) cases of $RS$'s it is possible to find a proper
scaling that secures  convergence.
 The disadvantage of approximation in this latter
 case is that the limiting process is no longer lies in $\cal N$.
 Some examples of scaling of logarithmic $RS$'s are presented e.g. in \cite{ABT}.
A key role here is played by Poisson-Dirichlet distribution on a
simplex as a limit of a scaled Ewens sampling formula. In
\cite{EG1}, it was proven a general result from which follows
(Corollary 3.1 there) that nondegenerate limiting distributions
are possible for convergent and logarithmic $RS$'s only. In the case
of  expansive $RS$'s the limiting distribution of a properly scaled
counting process is known to be a  curve called  a limit shape of
a random Young diagram (for the history of limit shapes see the
recent papers \cite{EG1}and \cite{Yaku}. Yakubovich(\cite{Yaku})
established general conditions for the existence of limit shapes
for scaled multiplicative measures that encompass known results on limit
shapes for regularly varying $RS$'s.

 \vskip.5cm
 {\bf Acknowledgement.} This research was supported
by New York Metropolitan Research Fund. A part of the present
paper was written during my visit to the School of Mathematical
Sciences, at Monash University. I am grateful to Fima Klebaner,
Aiden Sudbury and Kais Hamza for their hospitality.
I would also like to express my gratitude to a referee for important critical remarks and suggestions.

\end{document}